\documentclass[12pt]{amsart}
\usepackage{amssymb}
\usepackage{amsmath}
\usepackage{mathtools}
\usepackage[all]{xy}
\usepackage{enumitem}
\usepackage{xcolor}

 \rm

\renewcommand{\(}{\left(}
\renewcommand{\)}{\right)}
\renewcommand{\[}{\left[}
\renewcommand{\]}{\right]}

\newcommand{\abs}[1]{\left\lvert#1\right\rvert}
\newcommand{\norm}[1]{\left\lVert#1\right\rVert}
\newcommand{\st}{\:|\:}

\newcommand{\C}{{\mathbb{C}}}
\newcommand{\R}{{\mathbb{R}}}

\newcommand{\e}{\varepsilon}
\renewcommand{\phi}{\varphi}
\newcommand{\p}{\rho}

\newcommand{\supp}{{\mathrm{supp}}}
\newcommand{\tr}{{\mathrm{tr}}}

\renewcommand{\H}{{\mathcal{H}}}
\renewcommand{\S}{{\mathcal{S}}}
\newcommand{\BH}{{\mathcal{B}}(\H)}

\theoremstyle{plain}
\newtheorem{thm}{Theorem}[section]
\newtheorem{lem}[thm]{Lemma}

\theoremstyle{definition}

\theoremstyle{remark}

\title{Linear independence of quantum translates}

\author{Mansi Mishra}
\address{Department of Mathematics, IISER Bhopal, Bhopal 462 066.}

\begin{document}

\begin{abstract}
If $A$ is in the $p$-Schatten class on $\R^n$, $1\leq p \leq \frac{4n}{2n-1}$,
then the quantum translates of $A$ are linearly independent. Moreover, there
exists a non-zero operator in the $p$-Schatten class on $\R^n$,
$p>\frac{4n}{2n-1}$ whose quantum translates are linearly dependent.
\end{abstract}


\maketitle
\thispagestyle{empty}

\section{Introduction}

Let $f\in L^p(\R^n)$, $1\leq p \leq \infty$.  If $x\in \R^n$,
the translate of $f$ by $x$ is the function $T_xf:\R^n \to \C$ given by
$$
T_xf(y)= f(x+y), \quad y\in \R^n.
$$

In \cite{edgar},  Edgar and Rosenblatt considered the set of
translates of a non-zero function in $L^p(\R^n)$ by distinct elements
of $\R^n$, and gave conditions under which this set is linearly
independent.
They proved that any non-zero function in $L^p(\R)$,$1 \leq p < \infty$,
has linearly independent translates.
They also proved that if $n \geq 2$ and $1\leq p < 2n/(n-1)$,
then a non-zero function in $L^p(\R^n)$ has linearly independent translates
(see \cite[Corollary 2.7]{edgar}).
Moreover, they proved that there exists a non-zero function in $L^p(\R^n)$,
$n \geq 2$, that has linearly dependent translates when $p > 2n/(n-1)$
(see \cite[Example 2.4]{edgar}).
Furthermore, in \cite{rosenblatt}, Rosenblatt proved that a non-zero function
in $L^p(\R^n)$, $p=2n/(n-1)$, has linearly independent translates.

In this paper, we will prove non-commutative analogs of these results.
Recall that the reduced Heisenberg group $G$ is the set of triples
\begin{equation*}
\{(x,y,z) | x,y \in \R^n, z \in \C, \abs{z} = 1\}
\end{equation*}
with multiplication defined by
\begin{equation*}
(x,y,z)(x',y',z')=\(x+x', y+y', zz'e^{\pi i (x \cdot y'-y \cdot x')}\).
\end{equation*}
According to the Stone-von Neumann Theorem \cite[Theorem 1.50]{folland},
there is a unique irreducible unitary representation $\rho$ of $G$ such that
\begin{equation*}
\rho(0,0,z)= zI.
\end{equation*}
The standard realization of this representation is on the Hilbert space
$\H=L^2(\R^n)$ by the action
\begin{equation}\label{E:schrodinger}
\(\rho(x,y,z)\phi\)(t)= ze^{\pi i (x\cdot y+2y\cdot t)}\phi(t+x).
\end{equation}
Let $S^{p}(\H)$ denote the $p$-Schatten class of $\H$.
There is an action of $\R^{2n}$ on $S^p(\H)$ called quantum translation
which is defined by
\begin{equation*}
(x_1,y_1) \cdot X = \rho (x_1,y_1,1) X {\rho (x_1,y_1,1)}^{-1} .
\end{equation*}
In \cite{thangavelu-pw} and \cite{werner}, the quantum translation of
an operator was defined, and
it was shown that if $X \in S^2(\H)$ is non-zero,
then the quantum translates of $X$ are linearly independent.
The motivation for quantum translation may be found in \cite{thangavelu-pw}.

The following theorem is the main result of this paper, which is proved
in Section \ref{S:proof}.
\begin{thm}\label{T:main-thm}
Let $A \in S^p(\H)$, $1\le p \le 4n/(2n-1)$, be a non-zero operator and let
$(x_1,y_1), \dots, (x_k,y_k)$ be distinct elements of $\R^{2n}$.  Then
$\{(x_1,y_1) \cdot A, \dots, (x_k,y_k) \cdot A\}$ is a linearly independent set.
\end{thm}

In Section \ref{S:example}, we give an explicit example of a non-zero operator
in $S^p(\H)$, $p>4n/(2n-1)$, whose quantum translates are linearly dependent.

\section{Weyl Transform}

Recall that if $X \in S^1(\H)$, the {\em Fourier-Wigner transform} of $X$ is
the function on $\R^{2n}$ defined by
\begin{equation*}
\alpha(X)(x,y)= \tr(\rho(-x,-y,1)X)
\end{equation*}
and if $f\in L^1(\R^{2n})$, the {\em Weyl transform} of $f$ is the
operator $W(f) \in \BH$ defined by
\begin{equation*}
W(f) = \int_{\R^n} \int_{\R^n} f(x,y) \p(x,y,1) \,dx \,dy,
\end{equation*}
and that these are mutual inverses (see \cite{folland}).

We now discuss the relation between the quantum translates and the Weyl
transform (see \cite{thangavelu-pw} for a detailed discussion).
Recall that the Weyl correspondence of a function $f$ on $\R^{2n}$ is defined
by
$$
\operatorname{Weyl}(f)= W(\mathcal{F}^{-1}(f)),
$$
i.e., the Weyl transform of the inverse Fourier transform of $f$. It can be
shown that
$$
\operatorname{Weyl}(T_{(x,y)}f)= (x,y)\cdot \operatorname{Weyl}(f),
\quad (x,y)\in \R^{2n}.
$$

The definition of the Weyl transform has been extended to tempered
distributions by several authors (see e.g., \cite{maillard,keyl});
which we briefly recall now.
Let $\S(\R^{2n})$ denote the Schwartz class functions on $\R^{2n}$, and
let $\S'(\R^{2n})$ denote the space of tempered distributions on $\S(\R^{2n})$.
If $k \in \S(\R^{2n})$, the integral operator $X \in \BH$ defined by
\begin{equation*}
(X\phi)(x) = \int_{\R^n} k(x,y) \phi(y) \, dy
\end{equation*}
is called a Schwartz operator on $\H$.  The class of Schwartz
operators on $\H$ is denoted by $\S(\H)$, and is a Fr\'echet space
with respect to an appropriate sequence of seminorms.
Let $\S'(\H)$ denote the topological dual of $\S(\H)$.
It is well known that the Weyl transform is a topological
isomorphism from $\S(\R^{2n})$ to $\S(\H)$.
If $T \in \S'(\R^{2n})$, we define its Weyl transform $W(T) \in \S'(\H)$ by
\begin{equation*}
W(T)(Y)=T(\alpha(Y)^\sim), \quad Y\in \S(\H),
\end{equation*}
where $\alpha(Y)^\sim(x,y)=\alpha(Y)(-x,-y)$.
The distributional Weyl transform is a topological isomorphism from
$\S'(\R^{2n})$ to $\S'(\H)$, and its inverse is the distributional
Fourier-Wigner transform, which is defined by
\begin{equation*}
\alpha(X)(\phi) = X(W(\phi^\sim)), \qquad X\in\S'(\H), \quad \phi\in\S(\R^{2n}).
\end{equation*}

The following result was proved in \cite{mansi-wtod}.
\begin{thm}\label{T:wtod}
Let $T$ be a compactly supported distribution on $\R^{2n}$. Let $\widehat{T}$
denote the Fourier transform of $T$.
For $1 \leq p \leq \infty$, $W(T) \in S^p(\H)$ if and only if
$\widehat{T} \in L^p(\R^{2n})$. Moreover, if $K$ is a compact set in $\R^{2n}$,
then there exists a constant $C_K$ such that
\begin{equation*}
C_K^{-1} \norm{\widehat{T}}_p \leq
\norm{W(T)}_{S^p} \le C_K \norm{\widehat{T}}_p,
\end{equation*}
whenever $\supp(T)\subseteq K$.
Furthermore, $W(T)$ is compact if and only if
$\widehat{T} \in \mathcal{C}_0(\R^{2n})$.
\end{thm}

It is well known that if $\mu$ is a smooth measure on a hypersurface in $\R^n$,
$n \geq 2$, whose Gaussian curvature is nonzero everywhere, then
\begin{equation*}
\abs{\widehat{\mu}(\xi)} \leq A \abs{\xi}^{(1-n)/2},
\end{equation*}
where $A$ is a constant independent of $\xi$, and hence
$\widehat{\mu} \in L^p(\R^n)$ if $p > 2n/(n-1)$
(see e.g. \cite[p348, Theorem 1]{stein} and \cite[Theorem 7.7.14]{hormander}).
The following is an immediate consequence of Theorem \ref{T:wtod}.

\begin{thm}\label{C:hypersurface}
Suppose $S$ is a compact connected smooth hypersurface in $\R^{2n}$,
whose Gaussian curvature is nonzero everywhere.  
Let $\mu$ be a smooth measure on $S$.
Then $W(\mu)$ is a compact operator.
Moreover, $W(\mu) \in S^p(\H)$ if and only if $p> 4n/(2n-1)$.
\end{thm}

\section{Proof of the main result}\label{S:proof}

We prove Theorem \ref{T:main-thm} in this section. We will first recall
the proof of Theorem \ref{T:main-thm} for the case $p \in[1,2]$.

It can be shown that for any $A\in S^2(\H)$
\begin{equation}\label{E:weyl2}
\alpha((x_1,y_1) \cdot A)(x,y) = e((x_1,y_1),(x,y))\alpha(A)(x,y),
\end{equation}
where $e: \R^{2n} \times \R^{2n} \to \C$ is the map
\begin{equation*}
e((x',y'),(x,y))=e^{2\pi i(x'\cdot y- y'\cdot x)}.
\end{equation*}

Let $(x_1,y_1), \dots, (x_k,y_k)$ be distinct elements of $\R^{2n}$, and
let $c_1,\dots, c_k$ be nonzero scalars.
Define the difference operator $D$ on $S^p(\H)$, $1\leq p\leq \infty$, by
\begin{equation*}
DA = c_1 (x_1,y_1) \cdot A + \dots + c_k (x_k,y_k) \cdot A,
\quad A\in S^p(\H).
\end{equation*}
The trigonometric polynomial 
\begin{equation*}
c_1 e\((x_1,y_1),(x,y)\) + \dots + c_k e\((x_k,y_k),(x,y)\)
\end{equation*}
is called the characteristic trigonometric polynomial of the
difference equation $DA = 0$.

\begin{lem}
Let $A \in S^p(\H)$, $1\le p \le 2$, be a non-zero operator and let
$(x_1,y_1), \dots, (x_k,y_k)$ be distinct elements of $\R^{2n}$.  Then
$\{(x_1,y_1) \cdot A, \dots, (x_k,y_k) \cdot A\}$ is a linearly independent set.
\end{lem}  

\begin{proof}
Since $S^p(\H) \subseteq S^2(\H)$ for $1\le p \le 2$, it suffices to prove
the result for $p=2$.

Let $c_1,\dots, c_k \in \C$ be such that
\begin{equation*}
DA= c_1 (x_1,y_1) \cdot A + \dots + c_k (x_k,y_k) \cdot A = 0.
\end{equation*}
By the Plancherel theorem for $W$, we may write $A=W(f)$ for some nonzero
$f\in L^2(\R^{2n})$. 
By Equation (\ref{E:weyl2}), 
\begin{equation*}
DA = W\((c_1 e((x_1,y_1),\cdot)+ \dots + c_k e((x_k,y_k),\cdot)) f\).
\end{equation*}
Therefore $(c_1 e((x_1,y_1),\cdot)+ \dots + c_k e((x_k,y_k),\cdot)) f=0$.
Since $f$ is a non-zero function in $L^2(\R^{2n})$, there exists $\e >0$
such that $\abs{f} > \e$ on a set of positive measure.
Therefore the analytic function
\begin{equation*}
c_1 e((x_1,y_1),(x,y))+ \dots + c_k e((x_k,y_k),(x,y)) = 0
\end{equation*}
on a set of positive measure and thus it is identically zero.
Since
$$
\{e((x_1,y_1),\cdot), \dots, e((x_k,y_k),\cdot)\}
$$
is linearly independent, it follows that $c_1 = \dots = c_k = 0$.
\end{proof}
Let $A \in S^p(\H)$, $1 \leq p \leq \infty$, and let $T_A \in \S'(\H)$
be the map
$$
T_A(X) =\tr(XA).
$$
The Fourier-Wigner transform $\alpha(A)$ is defined in distributional sense,
i.e., $\alpha(A)= \alpha(T_A)$.
Let $\psi \in \S(\R^{2n})$.
\begin{equation*}
\begin{split}
&
\alpha((x_1,y_1) \cdot A)(\psi)\\
& =
T_{(x_1,y_1) \cdot A} (W({\psi}^\sim))\\
& =
\tr(W({\psi}^\sim)(x_1,y_1) \cdot A)\\
& =
\tr(W({\psi}^\sim)\rho (x_1,y_1,1) A {\rho (x_1,y_1,1)}^{-1})\\
& =
\tr\(\iint \psi(-x,-y) \rho(x,y,1) \rho (x_1,y_1,1) A
{\rho (x_1,y_1,1)}^{-1}\,dx \,dy\)\\
& =
\tr\(\iint e((x,y),(x_1,y_1)) \psi(-x,-y)\rho (x_1,y_1,1) \rho(x,y,1) A
{\rho (x_1,y_1,1)}^{-1}\,dx \,dy\)\\
& =
\tr\(\iint e((x,y),(x_1,y_1)) \psi(-x,-y) \rho(x,y,1) A \,dx \,dy\)\\
& =
\tr\(W(e((x_1,y_1), \cdot)\psi)^\sim)A\)\\
& =
\alpha(T_A)(e((x_1,y_1),\cdot) \psi)\\
& =
e((x_1,y_1),\cdot) \alpha(A)(\psi).
\end{split}
\end{equation*}
Therefore 
\begin{equation*}
\alpha(DA)=\(c_1 e((x_1,y_1),\cdot)+ \dots + c_k e((x_k,y_k),\cdot)\) \alpha(A).
\end{equation*}

Recall that the {\em symplectic Fourier transform} of a function
$f\in L^1(\R^{2n})$ is the function on $\R^{2n}$ defined by
\begin{equation*}
\breve{f}(\xi, \eta)=
\int_{\R^n}\int_{\R^n}  f(x,y) e^{2\pi i(\xi \cdot y - \eta \cdot x)} \,dx\,dy.
\end{equation*}
More generally (see \cite{maillard}), if $T$ is a tempered distribution
on $\R^{2n}$, the symplectic Fourier transform of $T$ is given by
\begin{equation*}
\breve{T}(\phi)= T((\breve\phi)^\sim),
\quad \phi\in \mathcal{S}(\R^{2n}).
\end{equation*}
Let $g(x,y)= e^{-\frac{\pi}{2}(\abs{x}^2+\abs{y}^2)}$, $(x,y) \in \R^{2n}$.
For $X \in \S'(\H)$, let $\beta(X)=g\alpha(X)$ and
$\breve{\beta}(X)$ denote the symplectic Fourier transform of $\beta(X)$.
It was shown in \cite{mansi-wtod} that $\breve{\beta}$ is a bounded linear
map from $S^p(\H)$ to $L^p(\R^{2n})$, $1 \leq p \leq \infty$. We will use
this fact in the proof of the main result.

Rosenblatt proved that the linear translates of a non-zero function in
$L^p(\R^n)$, $p= 2n/(n-1)$, are linearly independent by using the
following two results (see \cite{rosenblatt}).

\begin{lem}\label{L:rosenblatt}
If $f\in L^p(\R^n)$,  where $2<p \leq 2n/(n-1)$ and $n>1$, then
$$
\int_{\{R< \abs{x}<2R\}}\abs{f(x)}^2\, dx= o(R) \quad \text{as $R \to \infty$}.
$$
\end{lem}

\begin{thm}\label{T:rosenblatt}
Let $u$ be a tempered distribution supported on a real analytic set $K$ of
codimension $k>0$, and assume that $\hat{u}\in L^2_{\operatorname{loc}}(\R^n)$.
Assume that for every analytic manifold $M$ contained in $A$ and every
$x_0\in M$, there is an open cone $\Lambda \subset \R^n$ which contains some
normal of $M$ at $x_0$, and such that on
$\Lambda_R=\{x\in \Lambda \st R<\abs{x}<2R\}$,
$$
\lim_{R \to \infty}R^{-k}\int_{\Lambda_R} \abs{\hat{u}(x)}^2\,dx=0.
$$
Then $u=0$.
\end{thm}

We will now use these results to prove that quantum translates of a non-zero
operator in $S^p(\H)$, where $2< p \leq 4n/(2n-1)$, are linearly independent.

Let $A\in S^p(\H)$, $2< p \leq 4n/(2n-1)$, be a non-zero operator such that
$DA=0$.
It follows that 
\begin{equation}\label{E:da0}
0= \alpha(DA) =(c_1 e((x_1,y_1),\cdot)+ \dots + c_k e((x_k,y_k),\cdot))\alpha(A).
\end{equation}
Let $Z_k$ denote the zero set of the characteristic trigonometric polynomial
$c_1 e((x_1,y_1),\cdot)+ \dots + c_k e((x_k,y_k),\cdot)$.

\textbf{Case 1:}
Assume that $Z_k$ is of positive measure.
Since the analytic function
$c_1 e((x_1,y_1),(x,y))+ \dots + c_k e((x_k,y_k),(x,y))$
vanishes on a set of positive measure, it follows that it is identically zero.
Since the set
$$
\{e((x_1,y_1),\cdot), \dots, e((x_k,y_k),\cdot)\}
$$
is linearly independent,
it follows that $c_1 = \dots = c_k = 0$.
Therefore the quantum translates of $A$ are linearly independent.

\textbf{Case 2:}
If $Z_k$ is of measure zero, then $Z_k$ is an analytic submanifold of
codimension one (see e.g., \cite{complex}).
By Equation \ref{E:da0}, $\alpha(A)$ is supported in $Z_k$, and hence
$\beta(A)$ is supported in $Z_k$.
Take $u=\beta(A)$ in Theorem \ref{T:rosenblatt}.
Since $A\in S^p(\H)$, it follows that
$\breve{u}=\breve{\beta}(A) \in L^p(\R^{2n})$.
Therefore $\hat{u}\in L^2_{\operatorname{loc}}(\R^{2n})$.
By Lemma \ref{L:rosenblatt}, $u$ satisfies the growth condition required in
Theorem \ref{T:rosenblatt} for every cone $\Lambda$. Therefore $u=0$, and so
$\alpha(A)=0$. It follows that $A=0$.

\section{Example}\label{S:example}

We will now find an explicit example of a non-zero operator $A \in S^p(\H)$,
$p>4n/(2n-1)$, whose quantum translates are linearly dependent.
The argument used here is similar to the one given by Edgar and
Rosenblatt in \cite[Example 2.4]{edgar}.
The same example was given in \cite{wtom} to prove that there exists a non-zero
operator in $S^p(\H)$, $p>n \geq 6$, whose quantum translates are linearly
dependent, and in \cite{mansi-wtas} to prove that there exists a non-zero
compact operator on $\H$ whose quantum translates are linearly dependent.

Consider the difference equation given by
\begin{equation}\label{eq:de}
2(2n-1) A =
\sum_{j=1}^n \[(e_j, \mathbf{0})\cdot A + (-e_j, \mathbf{0})\cdot A
+ (\mathbf{0}, e_j)\cdot A +(\mathbf{0}, -e_j)\cdot A\],
\end{equation}
where $\mathbf{0}=(0, \dots, 0) \in \R^n$ and
$\{e_1, \dots, e_n\}$ is the standard basis of $\R^n$.

Then the characteristic trigonometric polynomial
\begin{equation*}
p(x_1, \dots, x_n, y_1, \dots, y_n) = 2(2n-1) - 2\sum_{j=1}^n \(\cos(2 \pi x_j) +
\cos(2\pi y_j)\)
\end{equation*}
has a zero set which is a disjoint union of compact $(2n-1)$-dimensional
surfaces of positive Gaussian curvature.  Let $S_{2n}$ be the connected
component of the zero set containing the points with all coordinates zero
except for one which is $\pm 1/4$.  Let $\sigma$ be the measure on $S_{2n}$
induced by the Lebesgue measure on $\R^{2n}$.
Then by Theorem \ref{C:hypersurface}, $W(\sigma) \in S^p(\H)$ for
$p> 4n/(2n-1)$.

It follows from equation (\ref{E:weyl2}) that $A=W(\sigma)$ is a
non-zero solution of the difference equation (\ref{eq:de}).
Hence, for $p> 4n/(2n-1)$, there exists an operator in $S^p(\H)$ whose
quantum translates are linearly dependent.

\bibliographystyle{amsplain}
\bibliography{v2-iqt}

\end{document}